\documentclass[a4paper,11pt]{article}
\usepackage{amsmath}
\usepackage{amssymb}
\usepackage{cite}
\usepackage{theorem}
{\theorembodyfont{\normalfont\slshape} \newtheorem{theorem}{Theorem}
\newtheorem{reftheorem}{Theorem}

\newtheorem{lemma}[theorem]{Lemma}

}
{\theorembodyfont{\normalfont\normalfont} }

\newcommand{\qed}{{$\qquad\square$\bigbreak}}
\newcommand{\integerset}{{\boldsymbol{Z}}}
\newcommand{\proof}[1]{%
\medbreak\noindent%
{\bfseries Proof{{\bfseries #1}}.\ \ }%
}

\newcommand{\DD}{\mathcal{D}}

\newcommand{\EEE}{\boldsymbol{E}}

\renewcommand{\baselinestretch}{1.2}
\begin{document}
\title{Precoloring Extension\\
Involving Pairs of Vertices of Small Distance}
\author{%
Chihoko Ojima\\
Akira Saito\footnote{%
Partially supported by Japan Society for the
Promotion of Science,
Grant-in-Aid
for Scientific Research (C),
22500018,
2012}
\\
Kazuki Sano\\
\ \ \\
Department of Information Science\\
Nihon University\\
Sakurajosui 3--25--40\\
Setagaya-Ku,
Tokyo 156--8550\\
JAPAN\\
E-mail : \texttt{asaito@chs.nihon-u.ac.jp}\ \ 
(Akira Saito)
}
\date{}
\maketitle
\renewcommand{\baselinestretch}{1.5}
\Large
\normalsize
\begin{abstract}
In this paper,
we consider coloring of graphs
under the assumption that some vertices are already colored.
Let $G$ be an $r$-colorable graph and let $P\subset V(G)$.
Albertson~[J.\ Combin.\ Theory Ser.~B
\textbf{73}
(1998),
189--194]
has proved that if every pair of vertices in $P$
have distance at least four,
then every $(r+1)$-coloring of $G[P]$ can be extended to
an $(r+1)$-coloring of $G$,
where $G[P]$ is the subgraph of $G$
induced by $P$.
In this paper,
we allow $P$ to have pairs of vertices of distance at most three,
and investigate how the number of such pairs affects the number
of colors we need to extend the coloring of $G[P]$.
We also study the effect of
pairs of vertices of distance at most two,
and extend the result by Albertson and Moore~[%
J.\ Combin.\ Theory Ser.~B
\textbf{77}
(1999)
83--95].
\end{abstract}
\bigbreak
\section{Introduction}
Graph coloring has a number of applications.
One example is a job scheduling problem.
In this problem,
each job is represented by a vertex,
and a pair of vertices are joined by an edge
if the corresponding jobs cannot be processed concurrently.
In this model,
an independent set represents a set of jobs which can be
preformed at the same time,
and if we assume that each job is processed in a unit time,
the chromatic number gives the minimum amount of time
in which we can finish all the jobs in a concurrent environment.
\par
In the real world,
however,
the job scheduling may not be tackled
from scratch.
In many cases,
the schedule of some jobs are already fixed and cannot be changed.
In graph coloring,
it corresponds to a situation in which some vertices are
already colored.
A precoloring extension is a problem to handle this situation.
In this problem,
a graph $G$,
a set of vertices $P\subset V(G)$ and
a coloring $d\colon P\to \integerset$ of
$G[P]$ are given,
where $G[P]$ is the subgraph of $G$ induced by $P$.
We call $d$ a \textit{precoloring}.
Our task is to find a coloring $V(G)\to\integerset$ of $G$ whose
restriction into $P$ coincides with $d$.
If $P$ is sufficiently sparse,
we may expect to extend $d$ to a coloring of $G$
with a few extra colors.
For the measure of sparseness,
Albertson~\cite{Albertson}
and
Albertson and Moore~\cite{AM}
have considered the minimum distance
between the vertices in $P$.
\begin{reftheorem}[\cite{Albertson}]\label{albertson}
Let $G$ be a graph with chromatic number at most $r$,
and
let $P\subset V(G)$.
Suppose every pair of distinct vertices in $P$ have distance
at least four.
Then every $(r+1)$-coloring of $P$ can be extended to an
$(r+1)$-coloring of $G$.
\end{reftheorem}
\begin{reftheorem}[\cite{AM}]\label{am}
Let $G$ be a graph with chromatic number at most $r$ and
let $P\subset V(G)$.
Suppose every pair of distinct vertices in $P$ have
distance at least three.
Then every $(r+1)$-coloring of $P$ can be extended to
a $\left\lceil\frac{3r+1}{2}\right\rceil$-coloring of $G$.
\end{reftheorem}
\par
These theorems give insight into
the relationship between the distance
of precolored vertices and the
number of colors necessary to extend
the precoloring to
a coloring of the whole graph.
\par
On the other hand,
again in the real world,
the assumptions of
Theorems~\ref{albertson}
and~\ref{am} may be idealistic.
For example,
while we have Theorem~\ref{albertson},
we may have to deal with a set $P$
of precolored vertices
which contains pairs of vertices of distance three.
In this case,
we might be forced to use more than $r+1$ colors to extend the
given precoloring.
But if the number of the pairs of distance three
is sufficiently small,
we expect that the number of additional colors
is also small.
Theorem~\ref{albertson}
does not answer this question.
\par
Motivated by this observation,
we investigate the situation in which
the set of precolored vertices contains pairs of distance
at most three and two,
and investigate how these pairs affect
the conclusion of Theorems~\ref{albertson}
and~\ref{am},
respectively
\par
For a graph $G$,
$P\subset V(G)$ and a positive integer $k$,
we define
$\DD_G(P, k)$ by
\[
\DD_G(P, k)=\bigl\{ \{x, y\}\subset P\colon
x\ne y \text{ and } d_G(x, y)\le k\bigr\},
\]
where $d_G(x, y)$ is the distance between $x$ and $y$ in $G$.
\par
In the next section,
we give an upper bound to the number of additional colors
to extend a given precoloring of $P$ to a coloring of $G$,
which is described in terms of
$|\DD_G(P, 3)|$.
In Section~3,
we give another bound,
which is described in terms of
$|\DD_G(P, 2)|$.
In Section~4,
we give some concluding remarks.
\par
We remark that this paper is neither the only nor the first
one to extend Theorems~\ref{albertson}
and~\ref{am}.
The problem of extending a precoloring
to the entire graph
has been studied in many papers.
We refer the readers who are interested in this problem
to~\cite{AH2, AH, AKW, AM2, AW, HT, HT2, HM, Tuza, TV, Voigt, Voigt2, Voigt3}.
In particular,
Albertson and Hutchinson~\cite{AH}
and Hutchinson and Moore~\cite{HM}
have considered the situation
in which the set of precolored vertices
induces a graph with several components,
and studied distance conditions among these components
that guarantee the extension
without using an additional color.
They have given best-possible results in many cases.
\par
For graph-theoretic notation and
definitions not explained in this paper,
we refer the reader to~\cite{CL}.
Let $G$ be a graph.
Then we denote by $\Delta(G)$ and $\chi(G)$
the maximum degree and the chromatic number of $G$,
respectively.
For $x\in V(G)$,
we denote the neighborhood of $x$ in $G$ by
$N_G(x)$.
In this paper,
we often deal with the closed neighborhood of $G$,
which is denoted and defined by
$N_G[x]=N_G(x)\cup\{x\}$.
If $A$,~$B\subset V(G)$ and $A\cap B=\emptyset$,
we define $E_G(A, B)$ by
$E_G(A, B)=\{ab\in E(G)\colon a\in A\text{ and }b\in B\}$.
\par
Let $P\subset V(G)$.
As we have already seen,
a coloring of $G[P]$ is called a precoloring of $P$ in $G$.
In this paper,
we always perceive a coloring of $G$ as
a mapping $f\colon V(G)\to\integerset$.
If $d\colon P\to\integerset$ is a precoloring of $P$ in $G$
and $f\colon V(G)\to\integerset$ is a coloring of $G$
with $f(v)=d(v)$ for every $v\in P$,
we say that $f$ \textit{extends\/} $d$.
For a positive integer $r$,
we denote the set $\{1, 2,\dots, r\}$ by $[r]$.
An $r$-coloring of $G$ is a coloring of $G$ which uses at most
$r$ colors.
In this paper,
an $r$-coloring is often perceived as a function
from $V(G)$ to $[r]$.
For $t\in [r]$,
$f^{-1}(t)$ is the set of vertices that receive the color~$t$.
We call it the \textit{color class\/}
of $V(G)$ with respect to the color~$t$.
\par
If $e=uv$ is an edge of a graph $G$,
we denote $\{u, v\}$ by $V(e)$.
Moreover,
for $F\subset E(G)$,
we write $V(F)$ for
$\bigcup_{e\in F}V(e)$.
A matching of $G$ is a set of independent edges in $G$.
Hence if $M$ is a matching,
then the order of $M$, denoted by $|M|$,
is the number of edges in $M$,
and $|V(M)|=2|M|$.
If $M$ is a maximum matching of $G$,
$|V(G)|-|V(M)|$ is called the deficiency of $G$.
Concerning the deficiency of a graph,
Berge's Formula is well-known.
We denote by $o(G)$ the number of components of odd order in $G$.
\begin{reftheorem}[Berge's Formula~\cite{Berge}]\label{berge}
For a graph $G$,
the deficiency of $G$ is given by $\max\{o(G-S)-|S|\colon
S\subset V(G)\}$.
\end{reftheorem}
\par
A matching $M$ in $G$ is called a perfect matching
if $V(M)=V(G)$,
and $M$ is called an almost perfect matching
if $|V(M)|=|V(G)|-1$.

\section{Pairs of Vertices of Distance Three}

In this section,
we investigate
the effect of the number of vertices
which are of distance at most three.
Theorem~\ref{albertson} states that
for a graph $G$ with
$\chi(G)\le r$ and $P\subset V(G)$
with
$\DD_G(P, 3)=\emptyset$,
every $(r+1)$-coloring of $P$
extends to an $(r+1)$-coloring of $G$.
If $\DD_G(P, 3)\ne\emptyset$,
we may need more than $r+1$ colors.
The purpose of this section is to prove that
for $t=|\DD_G(P,3)|$,
$r+O(\sqrt{t})$ colors suffice.

\begin{theorem}\label{main_distance_three}
Let $k$ be a positive integer.
Let $G$ be a graph with
$\chi(G)\le r$ and let $P\subset V(G)$.
Suppose
$|\DD_G(P, 3)| \le\frac{1}{2}k(k+1)$.
Then for each precoloring $d\colon P\to [r+1]$ in $G$,
there exists a coloring $f\colon V(G)\to [r+k]$ with
$f(u)=d(u)$ for each $u\in P$.
\end{theorem}

\par
We prove several lemmas to give a proof to
Theorem~\ref{main_distance_three}.
The following lemma has already been proved in~\cite{Albertson}.
But we tailor its statement so that it fits
the subsequent arguments.
For the completeness of the paper,
we give a proof to it.
For two colorings $f$,~g of a graph $G$,
we define $X(f, g)$ by
$X(f, g)=\{v\in V(G)\colon f(v)\ne g(v)\}$.

\begin{lemma}[\cite{Albertson}]\label{extension}
Let $G$ be a graph with an $r$-coloring
$c\colon V(G)\to [r]$.
Let $P$ be a set of vertices of $G$
with $\DD_G(P, 3)=\emptyset$.
Then for every precoloring
$d\colon P\to [r+1]$,
there exists an $(r+1)$-coloring $f\colon V(G)\to [r+1]$
of $G$ such that
\begin{enumerate}
\item
$f(x)=d(x)$
for every $x\in P$,
and
\item
for each $v\in X(c, f)$,
there exists a unique vertex $x\in P$ such that
$v\in N_G[x]$.
Moreover,
if $v\ne x$,
then $c(v)=d(x)$.
\end{enumerate}
\end{lemma}
\proof{}
For each $x\in P$,
if $d(x)\ne c(x)$,
then give the color $r+1$ to all the vertices $v$ in
$N_G(x)$ with $c(v)=d(x)$ and
then assign $d(x)$ to $x$.
Since $\DD_G(P, 3)=\emptyset$,
no two vertices receiving the color $r+1$ are adjacent.
Hence this gives a proper $(r+1)$-coloring $f$ of $G$.
By the construction,
$f$ satisfies both~(1)
and~(2).
\qed

\par
Let $k$ be a positive integer.
Then a coloring $c$ of a graph $G$ is said to
be an
\textit{almost $k$-coloring\/}
if
\begin{enumerate}
\item
$c$ is a $k$-coloring of $G$,
or
\item
$c$ is a $(k+1)$-coloring such that at least
one color class is a singleton set.
\end{enumerate}

\begin{lemma}\label{almost_r_coloring}
For a positive integer $k$,
a graph $H$ with $|E(H)|\le \frac{1}{2}k(k+1)$
has an almost $k$-coloring.
\end{lemma}
\proof{}
We proceed by induction on $k$.
If $k=1$,
then $|E(H)|\le 1$,
and it is easy to see that $H$ has
an almost $1$-coloring.
Suppose $k\ge 2$.
If $\chi(H)\le k$,
then a $k$-coloring of $H$ is an almost $k$-coloring of $H$.
Thus,
we may assume $\chi(H)\ge k+1$.
Then by Brook's Theorem,
we have $\Delta(H)\ge k$.
Let $x$ be a vertex of maximum degree,
and let $H'=H-x$.
Then $|E(H)|\le \frac{1}{2}k(k+1)-k=\frac{1}{2}k(k-1)$.
By the induction hypothesis,
$H$ has an almost $(k-1)$-coloring.
Then by assigning a new color to $x$,
we obtain an almost $k$-coloring of $H$.
\qed

\begin{lemma}\label{extension_2}
Let $G$ be a graph and let $P\subset V(G)$.
Suppose $P$ has a partition
$\{P_0, P_1,\dots, P_k\}$ with
$\DD_G(P_i, 3)=\emptyset$
for each $i$
with $0\le i\le k$.
Let $d\colon P\to [r+1]$ be
a precoloring of $P$ in $G$.
Suppose $G$ has an $r$-coloring $f$ with
$f(v)=d(v)$ for each $v\in P_0$.
Then for each $t$ with $0\le t\le k$,
there exists an $(r+t)$-coloring $f_t$ of $G$ with
$f_t(v)=d(v)$ for each $v\in \bigcup_{i=0}^t P_i$.
In particular,
$f_k$ is an $(r+k)$-coloring of $G$ with
$f_k(v)=d(v)$ for each $v\in P$.
\end{lemma}
\proof{}
We proceed by induction on $t$.
If $t=0$,
the lemma trivially follows with $f_0=f$.
Suppose $t\ge 1$.
By the induction hypothesis,
there exists an $(r+t-1)$-coloring $f_{t-1}$ of $G$
with $f_{t-1}(v)=d(v)$ for each $v\in\bigcup_{i=0}^{t-1}P_i$.
Without loss of generality,
we may assume
$f_{t-1}\colon V(G)\to [r+t-1]$.
By Lemma~\ref{extension},
there exists a coloring $f_t\colon V(G)\to [r+t]$ of $G$ such that
\begin{enumerate}
\item
$f_t(x)=d(x)$ for each $x\in P_t$,
and
\item
for each $v\in X(f_{t-1}, f_t)$,
there exists a unique vertex $x$ in $P_t$ with $v\in N_G[x]$.
Moreover,
if $v\ne x$,
then $f_{t-1}(v)=d(x)$.
\end{enumerate}
\par
Suppose $f_t(v)\ne d(v)$ for some $v\in\bigcup_{i=0}^t P_i$.
By~(1),
$v\in \bigcup_{i=0}^{t-1}P_i$ and hence
$f_{t-1}(v)=d(v)$.
This implies
$f_t(v)\ne f_{t-1}(v)$ and hence $v\in X(f_{t-1}, f_t)$.
Then there exists a unique vertex
$x\in P_t$ with $v\in N_G[x]$.
Since $v\notin P_t$,
we have $v\ne x$ and $f_{t-1}(v)=d(x)$
by~(2).
Now we have
$d(v)=f_{t-1}(v)=d(x)$,
which contradicts the assumption that $d$
is a precoloring of $P$.
Therefore,
we have
$f_t(v)=d(v)$ for each
$v\in\bigcup_{i=0}^t P_i$.
\qed
\proof{ of Theorem~\ref{main_distance_three}}
Define an auxiliary graph $H$ by
$V(H)=P$ and
$E(H)=\bigl\{uv\colon \{u, v\}\in \DD_G(P, 3)\bigr\}$.
Then $|E(H)|\le\frac{1}{2}k(k+1)$.
By Lemma~\ref{almost_r_coloring},
$H$ has an almost $k$-coloring.
Let $f_0\colon V(H)\to \{0\}\cup [k]$ be
an almost $r$-coloring of $H$ with
$|f_0^{-1}(0)|\le 1$.
Let $P_i=f_0^{-1}(i)$
($0\le i\le k$).
Since $\chi(G)\le r$,
there exists a coloring $c\colon V(G)\to [r]$ of $G$.
If $P_0\ne \emptyset$,
then let $P_0=\{u_0\}$ and
by taking an appropriate permutation of colors,
we may assume $c(u_0)=d(u_0)$.
Then by Lemma~\ref{extension_2},
$G$ has an $(r+k)$-coloring $f$ with
$f(v)=d(v)$ for each $v\in P$.
\qed

\section{Pairs of Vertices of Distance Two}

In this section,
we consider the effect of pairs of vertices
of distance two.
As we have seen in the introduction,
under the assumption of
$\DD_G(P, 2)=\emptyset$,
Albertson and Moore~\cite{AM} have proved that
an $(r+1)$-coloring of $P$ can be extended
to a $\left\lceil\frac{3r+1}{2}\right\rceil$-coloring of $G$.
They have also given
infinitely many examples of
$(G, P, d)$ for each positive integer $r$
such that
(1)~$G$ is a graph with $\chi(G)\le r$,
(2)~$P\subset V(G)$ with
$\DD_G(P, 2)=\emptyset$,
(3)~$d\colon P\to [r+1]$ is a precoloring of $P$ in $G$,
and
(4)~every coloring of $G$ that extends $d$
uses at least $\left\lceil\frac{3r+1}{2}\right\rceil$-colors.
Therefore,
Theorem~\ref{am} is best possible in this sense.
On the other hand,
the assumption on
$|\DD_G(P, 2)|$ has room for relaxation.
We prove that even if
$\DD_G(P, 2)\ne\emptyset$,
we can still extend an $(r+1)$-coloring of $P$
to a $\left\lceil\frac{3r+1}{2}\right\rceil$-coloring of $G$
as long as
$|\DD_G(P, 2)|$ is sufficiently small.
We also consider the case
in which $P$ is colored in more than $r+1$ colors.
\par
If $k\le r$,
we prove the following theorem.
\begin{theorem}\label{main_two_small_k}
Let $k$ and $r$ be positive integers with
$r\ge 2$ and $k\le r$,
and let $G$ be a graph with $\chi(G)\le r$.
Let $P\subset V(G)$
and let $d \colon P\to[r+k]$ be a precoloring of $P$ in $G$.
\begin{enumerate}
\item
If $r+k\equiv 0\pmod{2}$ and
$|\DD_G(P, 2)| < 2(r+k-1)$,
then $d$ can be extended to a
$\frac{3r+k}{2}$-coloring of $G$.
\item
If $r+k\equiv 1\pmod{2}$,
$r+k\ge 13$ and
$|\DD_G(P, 2)|< 3(r+k-1)$,
then $d$ can be extended to
a $\frac{3r+k+1}{2}$-coloring of $G$.
\end{enumerate}
\end{theorem}
We also prove that the bound of $|\DD_G(P,2)|$
in the above theorem is best-possible
in the following sense.
\begin{theorem}\label{sharpness_even}
For every pair of integers with $r$ and $k$ with
$r\ge 2$,
$r+k\equiv 0\pmod{2}$ and $k\le r$,
there exist infinitely many triples
$(G, P, d)$ such that
\begin{enumerate}
\item
$G$ is a graph of $\chi(G)\le r$,
\item
$P\subset V(G)$ and $|\DD_G(P, 2)| = 2(r+k-1)$,
\item
$d \colon P\to [r+k]$ is a precoloring of $P$ in $G$,
and
\item
$d$ cannot be extended to a $\frac{3r+k}{2}$-coloring of $G$.
\end{enumerate}
\end{theorem}
\begin{theorem}\label{sharpness_odd}
For every pair of integers $r$ and $k$ with
$r\ge 2$,
$r+k\equiv 1\pmod{2}$ and $k\le r$,
there exist infinitely many triples
$(G, P, d)$ such that
\begin{enumerate}
\item
$G$ is a graph with $\chi(G)\le r$,
\item
$P\subset V(G)$ and $|\DD_G(P, 2)|=3(r+k-1)$,
\item
$d \colon P\to [r+k]$ is a precoloring of $P$ in $G$,
and
\item
$d$ cannot be extended to a $\frac{3r+k+1}{2}$-coloring
of $G$.
\end{enumerate}
\end{theorem}
\par
In the range of $k > r$,
the situation changes.
We no longer need an additional color to extend
a precoloring of $P$.
Moreover,
different bounds of $|\DD_G(P, 2)|$
from those in Theorem~\ref{main_two_small_k}
appear.
\begin{theorem}\label{main_two_large_k}
Let $r$ and $k$ be positive integers with
$k > r\ge 2$,
and let $G$ be a graph with $\chi(G)\le r$.
Let $P\subset V(G)$ and let
$d \colon P\to [r+k]$ be a precoloring of $P$ in $G$.
\begin{enumerate}
\item
If $r < k \le\frac{3r-7}{2}$ and
$|\DD_G(P, 2)|<\min\left\{\frac{1}{2}(k+3r-4)(k-r+3),
(k-r+2)(k+r-1)\right\}$,
then $d$ can be extended to an $(r+k)$-coloring of $G$.
\item
If $k >\frac{3r-7}{2}$ and
$|\DD_G(P, 2)| < \min\left\{\frac{1}{2}(k+1)(k+2),
(k-r+2)(k+r-1)\right\}$,
then $d$ can be extended to an
$(r+k)$-coloring of $G$.
\end{enumerate}
\end{theorem}
\par
Note that the bounds of $|\DD_G(P, 2)|$ in Theorem~\ref{main_two_small_k}
are linear functions of $k$ and $r$ and they are sharp for $k \le r$,
while the bounds in the above theorem are quadratic functions
of $k$.
\par
In order to prove Theorems~\ref{main_two_small_k}
and~\ref{main_two_large_k},
we use matchings among colors.
For this purpose,
we introduce several definitions.
\par
For a set $X$,
a triple $\boldsymbol{E}=(E_0, E_1, E_2)$
is said to be
an \textit{ordered partition\/} of $X$ if
$E_0\cap E_1=E_0\cap E_2=E_1\cap E_2=\emptyset$
and $X=E_0\cup E_1\cup E_2$.
Note that in this definition,
possibly $E_i=\emptyset$ for some $i$,
$1\le i\le 3$.
Let $G$ be a graph,
and let $\boldsymbol{E}=(E_0, E_1, E_2)$
be an ordered partition of $E(G)$.
Let $M$ be a matching.
Then
$M$ is said to be a \textit{good matching\/}
for $\boldsymbol{E}$
if $M\cap E_2=\emptyset$
and $|M\cap E_1|\le 1$.
\par
Let $r$ and $k$ be positive integers.
Let $G$ be a graph with $\chi(G)\le r$ and let
$P\subset V(G)$.
Suppose a precoloring $d\colon P\to [r+k]$ in $G$ is given.
Now for a pair of distinct colors $\{i, j\}$,
we count how many times it appears as
$\{d(x), d(y)\}$ with
$\{x, y\}\in \DD_G(P, 2)$.
Let $H$ be the complete graph
with $V(H)=[r+k]$.
For each edge $e=ij$ of $H$,
we define $\varphi(e)$ by
$\varphi(e)=\bigl|\bigl\{ \{x, y\}\in \DD_G(P, 2)
\colon \{d(x), d(y)\}=\{i, j\}\bigr\}\bigr|$.
Then define $E_0^d$,
$E_1^d$,
$E_2^d$ and
$\boldsymbol{E}^d$ by
\begin{align*}
E_0^d &= \{e\in E(H)\colon \varphi(e)=0\},\\
E_1^d &= \{e\in E(H)\colon \varphi(e)=1\},\\
E_2^d &= \{e\in E(H)\colon \varphi(e)\ge2\}.
\quad\text{and}\\
\boldsymbol{E}^d &= (E_0^d, E_1^d, E_2^d).
\end{align*}
Clearly,
$\EEE^d$ is an ordered partition of $E(H)$.
\par
We make the following observation.
Though it is easy,
it is an important step in the proofs of Theorems~\ref{main_two_small_k}
and~\ref{main_two_large_k}.
\begin{lemma}\label{good}
Let $G$,
$P$,
$d$,
$H$,
$\varphi$ and $\EEE^d$
be as above.
Then
$|E_1^d|+2|E_2^d|\le|\DD_G(P, 2)|$.
\end{lemma}
\proof{}
By the definition of $\varphi$ and
$\EEE^d$,
\[
|\DD_G(P, 2)| = \sum_{e\in E(H)} \varphi(e)
=\sum_{k\ge 0}k|\{ e\in E(H)\colon \varphi(e)=k \}|
\ge |E_1^d|+2|E_2^d|.
\text{\qed}
\]
\par
We next prove a lemma which shows a relationship
between color extensions
and good matchings.
\begin{lemma}\label{good_matching}
Let $r$ and $k$ be positive integers,
and let $G$ be a graph with $\chi(G)\le r$.
Let $P\subset V(G)$ and let $d \colon P\to [r+k]$
be a precoloring of $P$ in $G$.
Let $H$ be the complete graph with $V(H)=[r+k]$.
\begin{enumerate}
\item
If $k\le r$ and $H$ has a good matching for
$\EEE^{d}$
of order $\left\lfloor\frac{1}{2}(r+k)\right\rfloor$,
then $d$ can be extended to
a $\left\lceil\frac{3r+k}{2}\right\rceil$-coloring of $G$.
\item
If $k > r$ and $H$ has a good matching
for $\EEE^{d}$ of order $r$,
then $d$ can be extended to
an $(r+k)$-coloring of $G$.
\end{enumerate}
\end{lemma}
\proof{}
Let $M=\left\{m_in_i\colon 1 \le i\le
\min\left\{r, \left\lfloor\frac{1}{2}(r+k)\right\rfloor
\right\}\right\}$ be a good matching for $\EEE^{d}$
in $H$.
We may assume $M\cap E_1^{d}\subset\{m_1n_1\}$,
and if $M\cap E_1^{d}=\{m_1n_1\}$,
then let $\{x_1, y_1\}$ be the unique pair in
$\DD_G(P, 2)$ with
$d(x_1)=m_1$ and $d(y_1)=n_1$.
\par
Since $\chi(G)\le r$,
there exists a coloring $g\colon V(G)\to [r]$
of $G$.
Let $D_i=g^{-1}(i)$
and $C_i=D_i-P$
($1\le i\le r$).
Without loss of generality,
we may assume $g(x_1)=1$.
Then $N_G(x_1)\cap C_1=\emptyset$.
\par
We first prove~(1).
We define a mapping $f\colon V(G)\to\integerset$ in the following way.
\begin{enumerate}
\item
If $v\in C_i$ with $1\le i\le \left\lfloor \frac{1}{2}(r+k)\right\rfloor$
and $N_G(v)\cap d^{-1}(m_i)=\emptyset$,
let $f(v)=m_i$.
\item
If $v\in C_i$ with
$1\le i\le \left\lfloor\frac{1}{2}(r+k)\right\rfloor$
and $N_G(v)\cap d^{-1}(m_i)\ne\emptyset$,
let $f(v)=n_i$.
\item
If $v\in C_i$ with
$\left\lfloor\frac{1}{2}(r+k)\right\rfloor+1\le i\le r$,
let $f(v)=\left\lceil\frac{1}{2}(r+k)\right\rceil+i$.
\item
If $v\in P$,
let $f(v)=d(v)$.
\end{enumerate}
Note that the total number of colors used by $f$
is at most $\left\lceil\frac{1}{2}(r+k)\right\rceil+r=
\left\lceil\frac{3r+k}{2}\right\rceil$.
\par
We prove that $f$ is a proper coloring of $G$.
Assume,
to the contrary,
that $f(u)=f(v)$ for some adjacent pair of vertices $u$ and $v$ in
$G$.
Since the restriction of $f$ into $P$ coincides with $d$
and $d$ is a precoloring,
$\{u, v\}\not\subset P$.
On the other hand,
since each $C_i$ is an independent set and
$f(C_i)\cap f(C_j)=\emptyset$ for each $i$ and $j$ with
$1\le i < j \le r$,
we have $\{u, v\}\not\subset \bigcup_{i=1}^r C_i$.
Thus,
$|\{u, v\}\cap P|=1$.
By symmetry,
we may assume $u\in P$ and $v\in C_i$ for some $i$,
$1\le i\le r$.
Then we have
$f(v)=f(u)=d(u)\le r+k$.
By the definition of $f$,
we have $i\le\left\lfloor\frac{1}{2}(r+k)\right\rfloor$
and $f(v)\in \{m_i, n_i\}$.
If $f(v)=m_i$,
then $d^{-1}(m_i)\cap N_G(v)=\emptyset$.
Hence $d(u)\ne m_i$.
Then $f(u)=d(u)\ne m_i=f(v)$.
This contradicts the assumption.
Thus,
$f(v)=n_i$.
This implies $d^{-1}(m_i)\cap N_G(v)\ne\emptyset$.
Let $w\in d^{-1}(m_i)\cap N_G(v)$.
Then $w\in P$ and $d(w)=m_i$.
\par
Since $d(u)=f(u)=f(v)=n_i$ and $d(w)=m_i\ne n_i$,
we have $u\ne w$.
Moreover,
since $v\in N_G(u)\cap N_G(w)$,
$d_G(u, w)\le 2$.
Therefore,
$\{u, w \}\in\DD_G(P, 2)$.
Then
$\varphi(m_in_i)\ge 1$.
Since $M$ is a good matching,
this is possible only if $i=1$,
$x_1=w$ and $y_1=u$.
However,
this implies $v\in N_G(x_1)\cap C_1$,
contradicting $N_G(x_1)\cap C_1=\emptyset$.
Therefore,
$f$ is a required coloring.
\par
Next,
we prove~(2).
Define a mapping $f\colon V(G)\to [r+k]$
in the following way.
\begin{enumerate}
\item
If $v\in C_i$,
$1\le i\le r$,
and $N_G(v)\cap d^{-1}(m_i)=\emptyset$,
let $f(v)=m_i$.
\item
If $v\in C_i$,
$1\le i\le r$,
and
$N_G(v)\cap d^{-1}(m_i)\ne\emptyset$,
let $f(v)=n_i$.
\item
If $v\in P$,
let $f(v)=d(v)$.
\end{enumerate}
Then the total number of colors used by $f$ is at most
$r+k$,
and the same argument as in the proof of~(1)
proves that $f$ is a proper coloring of $G$
\qed
\par
By the above lemma,
our main concern
is to find a good matching 
of an appropriate order
in the complete graph on a set of colors.
\begin{lemma}\label{even_and_small}
Let $r$ and $k$ be positive integers with
$r+k\equiv 0\pmod{2}$
and $k\le r$,
and let $H$ be the complete graph with
$V(H)=[r+k]$ and
let $\EEE=(E_0, E_1, E_2)$ be an ordered partition
of $E(H)$.
If
$|E_1|+2|E_2| < 2(r+k-1)$,
then $H$ contains a perfect matching which is also a good matching
for $\EEE$.
\end{lemma}
\proof{}
Since $H$ is a complete graph of even order,
$E(H)$ can be decomposed into $r+k-1$ perfect matchings
$M_1,\dots, M_{r+k-1}$.
Assume none of them is a good matching.
Then for each $i$,
$1\le i\le r+k-1$,
either $|M_i\cap E_1|\ge 2$ or
$|M_i\cap E_2|\ge 1$ holds.
Then in either case,
we have
$|M_i\cap E_1|+2|M_i\cap E_2|\ge 2$,
and
\[
\begin{split}
|E_1|+2|E_2| &= \sum_{i=1}^{r+k-1}|E_1\cap M_i|
+2\sum_{i=1}^{r+k-1}|E_2\cap M_i|\\
&=\sum_{i=1}^{r+k-1}(|M_i\cap E_1|+2|M_i\cap E_2|)
\ge 2(r+k-1).
\end{split}
\]
This contradicts the assumption,
and hence at least one of
$M_1,\dots, M_{r+k-1}$
is a good matching for $\EEE$.
\qed
By combining Lemmas~\ref{good},
\ref{good_matching}
and~\ref{even_and_small},
we obtain a proof of Theorem~\ref{main_two_small_k}~(1).
\par
Next,
we prove Theorem~\ref{main_two_small_k}~(2).
We further prove several lemmas.
\begin{lemma}\label{one_more}
Let $H$ be a complete graph,
and let $\EEE=(E_0, E_1, E_2)$
be an ordered partition of $E(H)$.
Let $M$ be a matching in $H$
with $M\subset E_0$ and
$|M| < \left\lfloor\frac{1}{2}|V(H)|\right\rfloor$,
and let $X=V(H)-V(M)$.
Then
\begin{enumerate}
\item
if $E(H[X])\cap (E_0\cup E_1)\ne\emptyset$
then $H$ contains a good matching for $\EEE$ of order
$|M|+1$,
and
\item
if $|E_H(V(M), X)\cap E_0| > |E_H(V(M), X)\cap E_2|$,
then $H$ contains a good matching for $\EEE$ of order
$|M|+1$.
\end{enumerate}
\end{lemma}
\proof{}
\noindent
(1)\ \ 
Let $e\in E(H[X])\cap (E_0\cup E_1)$.
Then $M\cup\{e\}$ is a good matching for $\EEE$.
\smallbreak\noindent
(2)\ \ Assume,
to the contrary,
that $H$ does not contain a good matching for $\EEE$
of order $|M|+1$.
Since
\[
\begin{split}
\sum_{e\in M}|E_H(V(e), X)\cap E_2)| &= |E_H(V(M), X)\cap E_2|\\
&< |E_H(V(M), X)\cap E_0|=\sum_{e\in M}|E_H(V(e), X)\cap E_0|,
\end{split}
\]
we have
$|E_H(V(e), X)\cap E_2| < |E_H(V(e), X)\cap E_0|$ for some $e\in M$.
Then since $|E_H(V(e), X)|=2|X|$,
we have $|E_H(V(e), X)\cap E_2| \le |X|-1$.
Let $e=uv$.
Since $E_H(V(e), X)\cap E_0\ne\emptyset$,
we may assume $ux_0\in E_0$
some $x_0\in X$.
If $vx\in E_0\cup E_1$ for some $x\in X-\{x_0\}$,
then
$(M-\{uv\})\cup\{ux_0, vx\}$ is a good matching for $\EEE$,
a contradiction.
Therefore,
$vx\in E_2$ for every $x\in X-\{x_0\}$.
Since
$|E_H(V(e), X)\cap E_2|\le |X|-1$,
This implies that
$E_H(V(e), X)\cap E_2=\bigl\{vx\colon x\in X-\{x_0\}\bigr\}$.
Then $|E_H(V(e), X)\cap E_0|\ge |X|$.
Since $|M|<\left\lfloor\frac{1}{2}|V(H)|\right\rfloor$,
$|X|\ge 2$.
Hence we can take a vertex $x_1\in X-\{x_0\}$.
Then $\{ux_1, vx_0\}\cap E_2=\emptyset$.
On the other hand,
since
$|E_H(V(e), X)\cap E_0|\ge |X|$,
$\{ux_1, vx_0\}\cap E_0\ne\emptyset$.
Then $(M-\{uv\})\cup\{ux_1, vx_0\}$ is a good matching
for $\EEE$.
This contradicts the assumption,
and the lemma follows.
\qed
\par
For a positive integer $n$,
we define an integer-valued function $h_n$
on $\integerset^+$ by
\[
h_n(t)=
\begin{cases}
\frac{1}{2}(t-1)(-t+2n) & \text{if $t < \frac{2}{5}(n+1)$}\\
(t-1)(2t-1) & \text{if $t \ge \frac{2}{5}(n+1)$}
\end{cases}
\]
\begin{lemma}\label{deficiency}
Let $n$ and $t$ be positive integers with
$t\le\left\lfloor \frac{1}{2}n\right\rfloor$.
Let $G$ be a graph of order $n$.
If $|E(G)| > h_n(t)$,
then $G$ contains a matching of order $t$.
\end{lemma}
\proof{}
Choose a graph $G_0$ of order $n$ without a matching of order $t$
so that $|E(G_0)|$ is as large as possible.
We prove $|E(G_0)|=h_n(t)$.
\par
By the assumption,
the deficiency of $G_0$
is at least $n-2(t-1)$.
By Berge's Formula,
$o(G_0-S)\ge |S|+n-2t+2$ for some $S\subset V(G)$.
Let $|S|=s$,
and let $C_1,\dots, C_k$ and $D_1,\dots, D_l$ be the odd
and even components of $G-S$,
respectively.
Since $t\le\left\lfloor \frac{1}{2}n\right\rfloor$,
$k\ge n-2t+2+s\ge 2$.
We may assume $|C_1|\ge |C_2|\ge\dots\ge |C_k|$.
\par
If $l\ge 1$,
then replace $C_1 \cup D_1$ with a complete graph of
order $|C_1|+|D_1|$.
Since $|C_1|+|D_1|$ is an odd number,
we have a new graph of order $n$ which has the same deficiency as that of
$G_0$ and contains more edges than $G_0$.
This contradicts the maximality of $|E(G_0)|$.
Therefore,
$l=0$,
and $G_0-S$ has no even components.
\par
If $k\ge s+n-2t+4$,
then $k\ge 3$.
Replace $C_1\cup C_2\cup C_3$ with a complete graph
of order $|C_1|+|C_2|+|C_3|$.
Let $G_1$ be the resulting graph.
Then $|V(G_1)|=n$,
$o(G_1-S)=k-2\ge s+n-2t+2$ and
$|E(G_1)| > |E(G_0)|$.
This again contradicts the maximality of $|E(G_0)|$.
Thus,
we have $k\le s+n-2t+3$.
However,
since $s+k\equiv n\pmod{2}$,
the equality does not hold,
and hence we have $k=s+n-2t+2$.
By the maximality of $|E(G_0)|$,
each of $C_1$,~$C_2,\dots, C_{s+n-2t+2}$ and
$S$ induces a complete graph.
Moreover,
there exists an edge between every vertex in $S$ and every vertex in
$\bigcup_{i=1}^{s+n-2t+2} C_i$.
\par
Assume $|C_2|\ge 3$ and replace $C_1$ and $C_2$ with complete graphs
of order $|C_1|+2$ and $|C_2|-2$.
Let $G_2$ be the resulting graph.
Then $G_2$ has the same deficiency as that of $G_0$,
and
\begin{multline*}
|E(G_2)| = |E(G_0)|-\frac{1}{2}|C_1|(|C_1|-1)
-\frac{1}{2}|C_2|(|C_2|-1)\\
+\frac{1}{2}(|C_1|+2)(|C_1|+1)+\frac{1}{2}(|C_2|-2)(|C_2|-3)
= |E(G_0)|+2(|C_1|-|C_2|+2).
\end{multline*}
Since $|C_1|\ge |C_2|$,
we have $|E(G_2)| > |E(G_0)|$.
This contradicts the maximality of $|E(G_0)|$.
Therefore,
we have
$|C_2|=|C_3|=\dots=|C_{s+n-2t+2}|=1$.
Then $|C_1|=n-(s+s+n-2t+2-1)=2t-2s-1$,
and hence
\[
\begin{split}
|E(G_0)| &= \frac{1}{2}(2t-2s-1)(2t-2s-2)
+\frac{1}{2}s(s-1)+s(n-s)\\
&=\frac{3}{2}s^2-\left(4t-\frac{5}{2}-n\right)s+(t-1)(2t-1).
\end{split}
\]
Let $f(s)=\frac{3}{2}s^2-\left(4t-\frac{5}{2}-n\right)s+(t-1)(2t-1)$.
Since $|C_1|\ge 1$,
$2t-2s-1\ge 1$,
and $s\le t-1$.
In the range of $0\le s\le t-1$,
$f(s)$ takes the maximum value at either $s=0$ or $s=t-1$.
Since $f(0)=(t-1)(2t-1)$ and
$f(t-1)=\frac{1}{2}(t-1)(-t+2n)$,
we have
$|E(G_0)|=\max\left\{(t-1)(2t-1), \frac{1}{2}(t-1)(-t+2n)\right\}=h_n(t)$.
\qed
\begin{lemma}\label{odd_and_small}
Let $r$ and $k$ be positive integers
with $r+k\equiv 1\pmod{2}$ and $k\le r$,
and let $H$ be the complete graph with $V(H)=[r+k]$.
Let $\EEE=(E_0, E_1, E_2)$ be an ordered partition of $E(H)$.
If $|E_1|+2|E_2|<\min\left\{3(r+k-1),\frac{1}{8}(r+k+3)(r+k+5)\right\}$,
then $H$ has an almost perfect matching
which is also a good matching for $\EEE$.
\end{lemma}
\proof{}
First,
we claim that there exists a matching $M$ with $M\subset E_0$ and
$|M|\ge \frac{1}{2}(r+k-3)$.
Assume the contrary.
Then by Lemma~\ref{deficiency},
we have
\[
|E_0|\le h_{r+k}\left(\tfrac{r+k-3}{2}\right)=
\begin{cases}
\frac{1}{2}(r+k-4)(r+k-5)&
\text{if $r+k\ge 19$}\\
\frac{3}{8}(r+k-5)(r+k+1) &
\text{if $r+k\le 17$.}
\end{cases}
\]
Then
\[
\begin{split}
|E_1|+2|E_2|\ge |E_1|+|E_2| & \ge\frac{1}{2}(r+k)(r+k-1)
-h_{r+k}\left(\tfrac{r+k-2}{2}\right)\\
& = \begin{cases}
4(r+k)-10 & \text{if $r+k\ge 19$}\\
\frac{1}{8}(r+k+3)(r+k+5) &
\text{if $r+k\le 17$.}
\end{cases}
\end{split}
\]
This immediately yields a contradiction
if $r+k\le 17$.
If $r+k\ge 19$,
we have
$4(r+k)-10\le |E_1|+2|E_2| < 3(r+k+1)$,
which yields $r+k < 13$,
again a contradiction.
Therefore,
the claim follows.
\par
If $|M|=\frac{1}{2}(r+k-1)$,
then $M$ is a required matching.
Therefore,
we may assume that $|M|=\frac{1}{2}(r+k-3)$
and $H$ does not have a good matching of order
$|M|+1$.
Let $X=V(G)-V(M)$.
Then $|X|=3$.
By Lemma~\ref{one_more},
we have
$E(H[X])\subset E_2$ and
$|E_H(V(M), X)\cap E_2|\ge|E_H(V(M), X)\cap E_0|$.
Then
\[
\begin{split}
|E_1|+2|E_2|&\ge |E_1\cap E_H(V(M), X)|
+2\bigl(|E_2\cap E_H(V(M), X)|+E(H[X])\bigr)\\
&= |E_1\cap E_H(V(M), X)|+2|E_2\cap E_H(V(M), X)|+2\cdot 3\\
&\ge |E_1\cap E_H(V(M), X)|+|E_2\cap E_H(V(M), X)|+
|E_0\cap E_H(V(M), X)|+6\\
&= |E_H(V(M), X)|+6=3(r+k-3)+6=3(r+k-1).
\end{split}
\]
This contradicts the assumption,
and the lemma follows.
\qed
\par
By combining Lemmas~\ref{good},
\ref{good_matching}
and~\ref{odd_and_small},
we obtain the following theorem,
which is slightly stronger than Theorem~\ref{main_two_small_k}~(2).
\begin{theorem}\label{main_extended}
Let $r$ and $k$ be positive integers with
$r\ge 2$,
$r+k\equiv 1\pmod{2}$ and $k < r$,
and let $G$ be a graph with $\chi(G)\le r$.
Let $P\subset V(G)$ and let $d\colon P\to [r+k]$
be a precoloring of $P$ in $G$.
If $|\DD_G(P, 2)|<
\min\left\{3(r+k-1), \frac{1}{8}(r+k+3)(r+k+5)\right\}$,
then $d$ can be extended to a
$\frac{3r+k+1}{2}$-coloring of $G$.
\end{theorem}
\par
Since
\[
\begin{split}
\min & \left\{ 3(r+k-1), \tfrac{1}{8}(r+k+3)(r+k+5)\right\}\\
& \ \ \ \ \ \ \ \ \ \ \ \ \ \ \ \ \ \ \ =
\begin{cases}
3(r+k-1) & \text{if $r+k\le 3$ or $r+k\ge 13$}\\
\frac{1}{8}(r+k+3)(r+k+5)
& \text{if $3 < r+k < 13$,}
\end{cases}
\end{split}
\]
$\frac{1}{8}(r+k+3)(r+k+5)$ appears as a bound
only if $r+k\in\{5, 7, 9, 11\}$.
The difference between
$\frac{1}{8}(r+k+3)(r+k+5)$ and
$3(r+k-1)$ in the range of $3\le r+k\le 13$ is shown
in the following table.
\begin{center}
\begin{tabular}{|c|cccccc|}
\hline
  & $3$ & $5$ & $7$ & $9$ & $11$ & $13$ \\
\hline
$3(r+k-1)$ & $6$ & $12$ & $18$ & $24$ & $30$ & $36$ \\
$\frac{1}{8}(r+k+3)(r+k+5)$ & $6$ & $10$ & $15$ & $21$ & $28$ & $36$\\
\hline
\end{tabular}
\end{center}
\medbreak
Next,
we prove Theorems~\ref{sharpness_even}
and~\ref{sharpness_odd}.
\proof{ of Theorem~\ref{sharpness_even}}
Let $q$ be an integer with $q\ge r+k$ and consider the following
construction.
\begin{enumerate}
\item
Construct a balanced complete $r$-partite graph $H$ with partite sets
$C_0$,~$C_1,\dots, C_{r-1}$,
where
$|C_0|=|C_1|=\dots=|C_{r-1}|=q$.
\item
Take a set of $2(r+k)$ vertices
$P=\{x_1,\dots, x_{r+k}, y_1,\dots, y_{r+k}\}$,
which is disjoint from $V(H)$.
\item
For each $i$,
$0\le i\le r-1$,
choose $r+k-1$ distinct vertices
$z_2^i$,~$z_3^i,\dots, z_{r+k}^i$ in $C_i$.
\item
Add edges
$\{x_1z_j^0, x_j z_j^0\colon 2\le j\le r+k\}
\cup\{y_1z_j^i, y_j z_j^i\colon 1\le i\le r-1, 2\le j\le r+k\}$.
\end{enumerate}
Let $G^q$ be the resulting graph.
Since $G^q$ is an $r$-partite graph with $C_0\cup\{y_1,\dots, y_{r+k}\}$,
$C_1\cup\{x_1,\dots, x_{r+k}\}$,
$C_2,\dots, C_{r-1}$,
we have $\chi(G^q)=r$.
By the construction,
we have
$\DD_{G^q}(P, 2)
=\bigl\{ \{x_1, x_j\}, \{y_1, y_j\}\colon 2\le j\le
r+k\bigr \}$
and hence
$|\DD_{G^q}(P, 2)|=2(r+k-1)$.
Define $d \colon P\to [r+k]$ by
$d(x_j)=d(y_j)=j$,
$1\le j\le r+k$.
\par
We prove that $d$ cannot be extended to
a $\frac{3r+k}{2}$-coloring of $G^q$.
Assume,
to the contrary,
that there exists a $\frac{3r+k}{2}$-coloring
$f\colon V(G^q)\to \left[\frac{3r+k}{2}\right]$
of $G^q$ such that the restriction of $f$ into $P$
coincides with $d$.
Choose $f$ so that
$\left|\bigcup_{i=0}^{r-1}f(C_i)\right|$ is as small as possible.
Since $H$ is a complete $r$-partite graph,
$f(C_i)\cap f(C_j)=\emptyset$ if $i\ne j$.
If $\{c_1, c_2\}\subset f(C_i)$
for some $c_1$ and $c_2$ with $c_1 >
r+k$ and $c_1\ne c_2$,
define $f'\colon V(G^q)\to\left[\frac{3r+k}{2}\right]$ by
\[
f'(v)=
\begin{cases}
c_1 & \text{if $v\in C_i$}\\
f(v) & \text{if $v\in V(G^q)-C_i$.}
\end{cases}
\]
Since $c_1 > r+k$,
$c_1\notin f(P)$.
Thus $f'$ is also a proper coloring of $G^q$ extending $d$.
Moreover,
$\bigcup_{i=0}^{r-1} f'(C_i)\subset \bigcup_{i=0}^{r-1}f(C_i)
-\{c_2\}$.
This contradicts the choice of $f$.
Therefore,
if $c\in f(C_i)$ for some $c > r+k$,
then $f(C_i)=\{c\}$.
Since $\frac{3r+k}{2}-(r+k)=\frac{r-k}{2}$,
at most $\frac{r-k}{2}$ partite sets of $H$ receive a color beyond $r+k$,
and hence
at least $\frac{r+k}{2}$ partite sets of $H$ are colored only
in colors in $[r+k]$.
Let $s=\frac{r+k}{2}$ and let
$C_{i_1}, C_{i_2},\dots, C_{i_s}$ be partite sets of $H$ with
$f(C_{i_j})\subset [r+k]$
($1\le j \le s$).
Since
$\bigl( N_{G^q}(x_l)\cup N_{G^q}(y_l)\bigr)\cap C_{i_j}\ne\emptyset$
for each $l$,
$1\le l\le r+k$,
we have $\left|f(C_{i_j})\right|\ge 2$
for each $j$,
$1\le j\le s$.
Therefore,
we have
$\left|\bigcup_{j=1}^s f(C_{i_j})\right|
=\sum_{j=1}^s |f(C_{i_j})|\ge 2s=r+k$.
This implies that
$\bigcup_{j=1}^s f(C_{i_j})=[r+k]$ and
$|f(C_{i_j})|=2$ for each $j$,
$1\le j\le s$.
Then $1\in f(C_{i_j})$ for some $j$.
Let $f(C_{i_j})=\{1, c\}$
($2\le c\le r+k$).
If $i_j\ne 0$,
then
$z_c^{i_j}\in N_{G^q}(y_1)\cap N_{G^q}(y_c)$ and hence
$f(z_c^{i_j})\cap \{1, c\}=\emptyset$.
On the other hand,
if $i_j=0$,
then
$z_c^0\in N_{G^q}(x_1)\cap N_{G^q}(x_c)$
and hence $f(z_c^0)\cap \{1, c\}=\emptyset$.
Therefore,
we have a contradiction in either case,
and hence $d$ cannot be extended to a
$\frac{3r+k}{2}$-coloring of $G^q$.
\qed
\proof{ of Theorem~\ref{sharpness_odd}}
Let $q$ be an integer with $r\ge 2$ and $q\ge r+k$.
We consider the following construction.
\begin{enumerate}
\item
Construct  a complete $r$-partite graph $H$ with
partite sets $C_0$,~$C_1,\dots, C_{r-1}$,
where
$|C_0|=|C_1|=\dots=|C_{r-1}|=q$.
\item
Take a set of $2(r+k)$ vertices
$P=\{x_1, x_2,\dots, x_{r+k}, y_1, y_2,\dots, y_{r+k} \}$
which is disjoint from
$\bigcup_{i=0}^{r-1}C_i$.
\item
Choose $r+k-2$ distinct vertices
$z_0^0$,~$z_4^0$,~$z_5^0,\dots, z_{r+k}^0$ from $C_0$.
\item
For every $i$ with $1\le i\le r-1$,
choose $r+k-3$ disjoint vertices
$z_4^i$,~$z_5^i,\dots, z_{r+k}^i$.
\item
Add edges
\begin{multline*}
\{x_1z_j^i, x_2z_j^i, x_3z_j^i, x_jz_j^i\colon
4\le j\le r+k,
1\le i\le r-1\}\\
\cup\{y_1z_0^0, y_2z_0^0, y_3z_0^0\}
\cup
\{y_jz_j^0\colon 4\le j\le r+k\}.
\end{multline*}
\end{enumerate}
Let $G^q$ be the resulting graph.
Define $d \colon P\to [r+k]$ by
$d(x_i)=d(y_i)=i$
($1\le i\le r+k$).
\par
Since $G^q$ is an $r$-partite graph
with partite sets
$C_0\cup\{x_1,\dots, x_{r+k}\}$,
$C_1\cup\{y_1,\dots, y_{r+k}\}$,
$C_2,\dots, C_{r-1}$,
we have $\chi(G^q)=r$.
By the definition of $P$,
\[
\begin{split}
\DD_{G^q}(P, 2)
&= \bigl\{ \{x_i, x_j\}\colon 1\le i\le 3, 4\le j\le r+k\bigr\}\\
&= \cup\bigl\{ \{x_1, x_2\}, \{x_1, x_3\}, \{x_2, x_3\}\bigr\}
\cup\bigl\{ \{y_1, y_2\}, \{y_1, y_3\}, \{y_2, y_3\}\bigr\}
\end{split}
\]
and hence
$|\DD_{G^q}(P, 2)|=3(r+k-3)+6=3(r+k-1)$.
\par
We prove that $G^q$ does not have a
$\frac{3r+k+1}{2}$-coloring which extends $d$.
Assume,
to the contrary,
that $G^q$ has a coloring $f\colon V(G^q)\to\left[\frac{3r+k+1}{2}\right]$
which extends $d$.
Choose $f$ so that
$\left|\bigcup_{i=0}^{r-1} f(C_i)\right|$ is as small as possible.
By the same argument as in the proof of Theorem~\ref{sharpness_even},
we have that
if $c\in f(C_i)$ for some $c > r+k$,
then $f(C_i)=\{c\}$.
\par
Since the number of colors
beyond $r+k$ is
$\frac{1}{2}(3r+k+1)-(r+k)=\frac{1}{2}(r-k+1)$,
these colors appear in at most
$\frac{1}{2}(r-k+1)$ partite sets of $H$.
Let $C_{i_1}$,~$C_{i_2},\dots, C_{i_{r'}}$ be the partite sets of
$H$ with
$f(C_{i_j})\subset [r+k]$
($1\le j\le r'$).
Then $r'\ge r-\frac{1}{2}(r-k+1)=\frac{1}{2}(r+k-1)$.
\par
By the construction,
$f(C_{i_s})\cap f(C_{i_t})=\emptyset$
for each $s$,~$t$ with
$1\le s < t \le r'$.
Moreover,
$\{y_1, y_2, y_3\}\subset N_{G^q}(z_0^0)$,
$y_j\in N_{G^q}(z_j^0)$
and
$\{x_1, x_2, x_3, x_j\}\subset N_{G^q}(z_j^i)$
for every $i$,~$j$
with
$1\le i\le r-1$ and $4\le j\le r+k$.
These imply $|f(C_{i_t})|\ge 2$ for each $t$,
$1\le t\le r'$.
Therefore,
\[
r+k \ge \left|\bigcup_{j=1}^{r'}f(C_{i_j})\right|
=\sum_{j=1}^{r'}|f(C_{i_j})|
\ge 2r' \ge r+k-1
\]
Since $r+k\equiv 1\pmod{2}$,
we have $r'=\frac{1}{2}(r+k-1)$.
\par
Suppose
$|f(C_{i_1})|=|f(C_{i_2})|=\dots=|f(C_{i_{r'}})|=2$.
Then
$\left|\bigcup_{j=1}^{r'} f(C_{i_j})\right|=2r'=r+k-1$.
Let $c_0$ be the unique color in
$[r+k]-\bigcup_{j=1}^{r'} f(C_{i_j})$.
By the symmetry of the colors $1$,~$2$,
and~$3$,
we may assume
$c_0\notin \{1, 2\}$
and $1\in f(C_{i_1})$.
Let $f(C_{i_1})=\{1, c_1\}$.
If $i_1\ne 0$,
consider $z_{c_1}^{i_1}\in C_{i_1}$.
Since
$\{x_1, x_{c_1}\}\subset N_{G^q}(z_{c_1}^{i_1})$,
$f(x_1)=d(x_1)=1$ and
$f(x_{c_1})=d(x_{c_1})=c_1$,
neither $1$ nor $c_1$ is allowed for
$f(z_{c_1}^{i_1})$,
a contradiction.
Therefore,
$i_1=0$ and $1\in f(C_0)$.
By the same argument,
we also have $2\in f(C_0)$,
and hence $f(C_0)=\{1, 2\}$.
However,
since $\{y_1, y_2\}\subset N_{G^q}(z_0^0)$,
$d(y_1)=1$ and
$d(y_2)=2$,
neither $1$ nor $2$ can be assigned to $z_0^0$.
This is again a contradiction.
Therefore,
$|f(C_{i_j})|\ge 3$ for some $j$,
$1\le j\le r'$.
Then we may assume
$|f(C_{i_1})|=3$
and
$|f(C_{i_j})|=2$ for each $j$ with
$2\le j\le r'$.
This implies
$\bigcup_{j=1}^{r'} f(C_{i_j})=[r+k]$.
\par
Assume $|f(C_{i_1})\cap \{1, 2, 3\}|\le 1$.
Then we may assume
$\{1, 2\}\cap f(C_{i_1})=\emptyset$.
Let $1\in f(C_{i_j})$,
$2\le j\le r'$,
and let $f(C_{i_j})=\{1, c_1\}$.
By the same argument as in the previous paragraph,
we have $i_j=0$ and $1\in f(C_0)$.
Similarly,
$2\in f(C_0)$.
Then $f(C_0)=\{1, 2\}$,
a contradiction.
Therefore,
$|f(C_{i_1})\cap \{1, 2, 3\}|\ge 2$.
Since
$\{x_1, x_2, x_3\}\subset\bigcap_{i=1}^{r-1}N_{G^q}(z_1^i)$ and
$\{y_1, y_2, y_3\}\subset N_{G^q}(z_0^0)$,
$f(C_{i_1})\ne\{1, 2, 3\}$.
Therefore,
we have
$f(C_{i_1})=\{1, 2, c\}$ for some $c$ with
$4\le c\le r+k$.
Since $\{x_1, x_2, x_c\}\subset\bigcap_{i=1}^{r-1}N_{G^q}(z_c^i)$,
$i_1\notin\{1, 2,\dots, r-1\}$
and hence $i_1=0$,
or
$f(C_0)=\{1, 2, c\}$.
This implies
$f(C_{i_j})=\{3, c'\}$ for some $j$,
$2\le j\le r'$,
and $c'\in [r+k]$.
However,
since $\{x_3, x_{c'}\}\subset \bigcap_{i=1}^{r-1}N_{G^q}(z_{c'}^i)$,
$i_j\notin \{1,\dots, r-1\}$ and hence $i_j=0$.
This means $f(C_0)=\{3, c'\}$.
This is a final contradiction,
and the theorem follows.
\qed
Now we prove Theorem~\ref{main_two_large_k}.
In view of Lemmas~\ref{good} and~\ref{good_matching},
it suffices to prove the following lemma.
\begin{lemma}
Let $r$ and $k$ be positive integers with $k > r$,
and let $H$ be the complete graph of order $r+k$.
Let $\EEE=(E_0, E_1, E_2)$ be an ordered partition of $E(H)$.
Then
\begin{enumerate}
\item
if $r < k \le \frac{3r-7}{2}$ and
$|E_1|+2|E_2| <\min\left\{\frac{1}{2}(k+3r-4)(k-r+3),
(k-r+2)(k+r-1)\right\}$,
then
$H$ has a good matching of order $r$,
and
\item
if $k> \frac{3r-7}{2}$ and
$|E_1|+2|E_2| <\min\left\{\frac{1}{2}(k+1)(k+2),
(k-r+2)(k+r-1)\right\}$,
then $H$ has a good matching of order $r$.
\end{enumerate}
\end{lemma}
\proof{}
Note
\[
h_{r+k}(r-1)=
\begin{cases}
(2r-3)(r-2) & \text{if $r < k \le \frac{3r-7}{2}$}\\
\frac{1}{2}(r-2)(r+2k+1) & \text{if $k > \frac{3r-7}{2}$.}
\end{cases}
\]
Let $\bar{h}(r, k)=\frac{1}{2}(r+k)(r+k-1)-h_{r+k}(r-1)$.
Then
\[
\bar{h}(r, k)=
\begin{cases}
\frac{1}{2}(k+3r-4)(k-r+3) & \text{if $r < k\le \frac{3r-7}{2}$}\\
\frac{1}{2}(k+1)(k+2) & \text{if $k > \frac{3r-7}{2}$.}
\end{cases}
\]
Hence we have $|E_1|+2|E_2| < \bar{h}(r, k)$
both in~(1) and~(2).
\par
Let $H_0=(V(H), E_0)$.
Since $|E_0|+|E_1|+|E_2|=\frac{1}{2}(r+k)(r+k-1)$,
we have
\[
\begin{split}
|E(H_0)| &=\frac{1}{2}(r+k)(r+k-1)-(|E_1|+|E_2|)\\
&\ge \frac{1}{2}(r+k)(r+k-1)-(|E_1|+2|E_2|)\\
&>\frac{1}{2}(r+k)(r+k-1)-\bar{h}(r, k)
=h_{r+k}(r-1),
\end{split}
\]
$H_0$ has a matching $M_0$ with $|M_0|\ge r-1$.
\par
Assume $H$ does not have a good matching for $\EEE$
of order $r$,
Then $|M_0|=r-1$.
Let $X=V(H)-V(M_0)$.
Then
$|X|=(r+k)-2(r-1)=k-r+2\ge 3$
since $k > r$.
By Lemma~\ref{one_more},
$E(H[X])\subset E_2$
and
$|E_H(V(M), X)\cap E_2|\ge |E_H(V(M), X)\cap E_0|$.
Therefore,
\[
\begin{split}
|E_1|+2|E_2| & \ge
|E_H(V(M), X)\cap E_1|+2\bigl(|E_H(V(M), X)\cap E_2|+|E(H[X])|\bigr)\\
&=|E_H(V(M), X)\cap E_1|+2|E_H(V(M), X)\cap E_2|
+2|E(H[X])|\\
&\ge |E_H(V(M), X)\cap E_1|+|E_H(V(M), X)\cap E_2|\\
&\ \ \ \ \ \ \ \ \ \ \ \ \ \ \ \ \ \ \ \ 
\ \ \ \ \ \ \ \ \ \ \ \ \ \ +|E_H(V(M), X)\cap E_0|+2|E(H[X])|\\
&=|E_H(V(M), X)|+2\cdot\frac{1}{2}(k-r+2)(k-r+1)\\
&=2(r-1)(k-r+2)+(k-r+2)(k-r+1)\\
&=(k-r+2)(k+r-1).
\end{split}
\]
This contradicts the assumption,
and the lemma follows.
\qed

\section{Conclusion}

In this paper,
we have studied extension of precoloring.
For a set of precolored vertices $P$ in a graph $G$,
we define $\DD_G(P, k)$ to be the set of pairs of vertices in $P$
whose distance is at most $k$.
We have investigated how $|\DD_G(P, 3)|$ and $|\DD_G(P, 2)|$
affect the bounds on the number of additional colors required
to extend the precoloring of $P$,
which have been given by Albertson~\cite{Albertson} and
Albertson and Moore~\cite{AM}.
We have proved that the sharpness of Theorem~\ref{main_two_small_k}.
However,
the sharpness of Theorems~\ref{main_distance_three}
and~\ref{main_two_large_k}
is unknown.

\end{document}